\documentclass[a4paper,10pt]{article}

\usepackage{amssymb,amsmath,graphics}

\newenvironment{dwd}{\par\noindent{\bf Proof.}}{\par\rightline{$\blacksquare$}}

\newtheorem{theo}{Theorem}
\newtheorem{prop}{Proposition}  
\newtheorem{coro}{Corollary}

\newtheorem{rema}{Remark}
\newtheorem{lema}{Lemma}
\newtheorem{defi}{Definition}
\def\be#1\ee{\begin{equation}#1\end{equation}}
\newcommand{\ba}{\begin{eqnarray} }
\newcommand{\ea}{\end{eqnarray} }
\def\bt#1\et{\begin{theo}#1\end{theo}}
\def\bl#1\el{\begin{lema}#1\end{lema}}
\def\bp#1\ep{\begin{prop}#1\end{prop}}
\def\bd#1\ed{\begin{defi}#1\end{defi}}

\def\ccB{{\cal B}}

\def\ccF{{\cal F}}

\def\ccX{{\cal X}}

\def\ra{\rightarrow}

\def\E{\mathbf{E}}
\def\P{\mathbf{P}}

\def\R{{\mathbb R}}

\def\ls{\leqslant}
\def\gs{\geqslant}

\begin{document}

\title{\bf Concentration via chaining method and its applications}
\author{{Witold Bednorz}
\footnote{{\bf Subject classification:} 60G15, 60G17}
\footnote{{\bf Keywords and phrases:} sample boundedness, Gaussian processes}
\footnote{Support: Polish Ministry of Science and Higher Education Iuventus Plus Grant no. IP 2011 000171}
\footnote{Institute of Mathematics, University of Warsaw, Banacha 2, 02-097 Warszawa, Poland}}

\maketitle
\begin{abstract}
In this paper we study the regularity of paths in terms of properties of admissible nets. We show the 
right concentration inequality above the modulus of continuity. Using the approach we prove the Bernstein type inequality for the empirical processes. 
Therefore we obtain the best form of concentration for processes studied recently by Mendelson and Paouris and Tomczak-Jaegerman. 
Results of this type are of importance in the compressed sensing theory. 
\end{abstract}

\section{Introduction}

In this paper we show how the chaining approach can  be used to establish results on the concentration of 
functionals of random variables. In particular we focus on the well known problem in the
theory of empirical processes. We show how to prove the right concentration inequality on the supremum of
centered sums of squares of independent random variables. 
The chaining approach has been recently used to investigate empirical processes \cite{MPT} and \cite{PaM}.
The main difficulty that appears in the results is to use a special chaining depending 
on the approximation level. We start our study form the simplest case of
a stochastic process with increments under control of one distance. The we turn to more involved case of two 
distances. With the results we turn to prove the application to the empirical processes proving the Bernstein type 
inequality. Finally we discuss how the result can be used in the compressed sensing. 

\section{One distance control of increments}

Let $(T,d)$ be a compact metric space. Let $X(t)$, $t\in T$ be a stochastic process defined on $(T,d)$. We aim to study path properties of $X(t)$, $t\in T$
under some increment conditions. The simplest setting in which the problem can be analyzed is when there is a single distance $d$ on $T$
and a single Young function $\psi:\R_{+}\ra \R_{+}$, i.e. convex, increasing and such that $\psi(0)=0$ and $\psi(1)=1$ which are used to
impose the following control on the increments
\be\label{erm1}
\E\psi(\frac{|X(t)-X(s)|}{d(s,t)})\ls 1,\;\;\mbox{for all}\;s,t\in T.
\ee 
Note that defining for all $s,t\in T$
\be\label{erm2}
\|X(t)-X(s)\|_{\psi}=\inf\{C>0:\;\E\psi(\frac{|X(t)-X(s)|}{C})\ls 1\}<\infty
\ee
one can use distance $d$ of the form $d(s,t)=\|X(t)-X(s)\|_{\psi}$. Obviously the condition (\ref{erm1}) implies some concentration
inequality for increments, i.e.
$$
\P(|X(t)-X(s)|\gs ud(s,t))\ls \E \psi(\frac{|X(t)-X(s)|}{d(s,t)})/\psi(u)\ls \frac{1}{\psi(u)}\;\;\mbox{for}\;u>0.
$$
Note that the requirement that $\psi$ is convex can be slightly relaxed to 
the condition that $\psi$ is continuous, increasing to infinity, $\psi(0)=0$ and $\psi(1)=1$. It usually concerns the case  when 
the convexity starts for large enough arguments.
\smallskip

\noindent
The simplest example of such an increment control is when $\psi$ is exponential, i.e.
we assume that $\psi(x)=\varphi_p(x)=2^{x^{p}}-1$, $p>0$. Then (\ref{erm1}) is equivalent to the
following concentration inequality
$$
\P(|X(t)-X(s)|>d(s,t)u)\ls A\exp(-Bu^p),\;\;\mbox{for all}\;u\gs 0,
$$
where $A,B$ - universal constants. In order to analyze the whole class of exponential functions $\psi$
we formulate the general exponential-type condition such functions should satisfy
\be\label{erm3} 
\psi^{-1}(xy)\ls K(\psi^{-1}(x)+\psi^{-1}(y))\;\;x,y\gs 0,
\ee
where $K<\infty$ is a universal constant. Note that $K=1$ for $\psi=\varphi_p$, $p\gs 1$, nevertheless the inequality holds
for all $\psi=\varphi_p$, $p>0$ but for $p<1$, $K$ depends on $p$.
\smallskip

\noindent
In the particular case of centered Gaussian processes $X(t)$, $t\in T$ one can apply  
$\psi(x)=\varphi_2(x)=2^{x^2}-1$ and $d(s,t)=\sqrt{8/(3\ln 2)}(\E|X(t)-X(s)|^2)^{\frac{1}{2}}$. 
The meaning of (\ref{erm1}) is that we acknowledge Gaussian-type concentration for increments, as
the best possible tool to understand path properties of the process.

\section{The concentration inequality for one distance }

We start our study towards concentration inequalities by the chaining approach
from the simplest setting of a compact $(T,d)$ and a single Young function $\psi$.
It should be stressed that if $T$ is not a completely bounded space in $d$ one can construct 
processes that satisfies (\ref{erm1}) and which are not sample bounded.   
\smallskip

\noindent
In this setting there exists a separable modification of $X(t)$, $t\in T$ which we refer to from now on.
Indeed (\ref{erm1}) implies the continuity in probability of the process
and since $(T,d)$ is compact we can define a separable modification of $X(t)$, $t\in T$ based on any separable dense subset $T_0\subset T$.
\smallskip

\noindent
Let $(T,d)$ be a compact metric space. 
Let $N_n=\psi(2^{n})$, for $n\gs 1$ and $N_0=1$. 
We say that a sequence of finite sets $(T_n)^{\infty}_{n=0}$ is admissible, if $|T_n|\ls N_{n}$, $T_{n}\subset T_{n+1}$
and $\bigcup^{\infty}_{n=0}T_n$ is dense in $T$.
\smallskip

\noindent
For each $n\gs 0$ define $\pi_n(t)\in T_n$ in a way that $d(t,T_n)=d(t,\pi_n(t))$.
The first question we study is what one can say on the
difference $|X(t)-X(\pi_m(t))|$ for a given $m\gs 0$.
For each $m\gs 0$ define
$$
\sigma_m(t)=\sum^{\infty}_{n=m}2^nd(t,T_n).
$$
Moreover let 
$$
Z=\sum^{\infty}_{n=1}\frac{1}{N_n^3}\sum_{t\in T_{n}}\psi(\frac{|X(t)-X(\pi_{n-1}(t))|}{d(t,\pi_{n-1}(t))}).
$$
Clearly $\E Z\ls 1$, since $\frac{N^2_n}{N_{n}^3}\ls \frac{1}{N_n}\ls \frac{1}{2^n}$, by the convexity of $\psi$. 
The following result is a stronger version of the usual concentration inequality obtained by the chaining approach (cf. \cite{Tal}).
\begin{prop}\label{pro1}
For each separable $X(t)$, $t\in T$ that satisfies (\ref{erm1}) for $\psi$ that verifies the exponential condition (\ref{erm3})
the following holds. For all $t\in T$ and $m\gs 0$ 
$$
|X(t)-X(\pi_m(t))|\ls A\sigma_{m}(t)+Bd(t,\pi_m(t))\psi^{-1}(Z),
$$
where $A=15K^2,B=4K^2$ are universal constants, $Z\gs 0$ and $\E Z\ls 1$.
\end{prop}
\begin{dwd}
Fix $t\in T$. We define $n_{i}=n_i(t)$, $i\gs 0$ that may depend on $t$ in the following way.
Let $n_0=m$ and for $i\gs 1$ let $n_i=\inf\{n>n_{i-1}:\;2d(t,T_n)< d(t,T_{n_{i-1}})\}$
We use the following chaining argument
\be\label{a2}
|X(t)-X(\pi_m(t))|\ls \sum^{\infty}_{i=1}|X(\pi_{n_i}(t))-X(\pi_{n_{i-1}}(t))|,
\ee
Consequently due to (\ref{erm3})
\begin{align*}
& |X(\pi_{n_i}(t))-X(\pi_{n_{i-1}}(t))|\ls\\
&\ls d(\pi_{n_{i}}(t),\pi_{n_{i-1}}(t))\psi^{-1}(\frac{N_{n_{i}}^3}{N_{n_{i}}^3}V_{n_i})=\\
&\ls  d(\pi_{n_{i}}(t),\pi_{n_{i-1}}(t))[3K^2 2^{n_{i}}+K^2\psi^{-1}(\frac{1}{N_{n_{i}}^3}V_{n_i})],
\end{align*}
where
$$
V_{n}=\sum_{u,v\in T_{n}}\psi(\frac{|X(u)-X(u)|}{d(u,v)})\;\;\mbox{for all}\;n\gs 0.
$$
By the definition of $n_i$, $i\gs 0$
$$
2d(t,T_{n_i-1})\gs d(t,T_{n_{i-1}}). 
$$
Therefore
\begin{align*}
& d(\pi_{n_{i}}(t),\pi_{n_{i-1}}(t))\ls d(t,\pi_{n_i}(t))+\\
&+d(t,\pi_{n_{i-1}}(t))\ls d(t,T_{n_i})+2d(t,T_{n_{i}-1})
\end{align*}
and hence
\begin{align*}
& 3K^2\sum^{\infty}_{i=1} 2^{n_i}d(\pi_{n_{i}}(t),\pi_{n_{i-1}}(t))\ls 
3K^2 \sum^{\infty}_{i=1} 2^{n_i}(d(t,T_{n_{i}})+2d(t,T_{n_{i}-1}))\ls\\
&\ls 3K^2(\sum^{\infty}_{n=0} 2^n d(t,T_n)+4\sum^{\infty}_{n=0} 2^n d(t,T_n))=
15K^2 \sigma_m(t).
\end{align*}
On the other hand 
\begin{align*}
&\sum^{\infty}_{i=1}d(\pi_{n_{i}}(t),\pi_{n_{i-1}}(t))\ls 2 \sum^{\infty}_{i=0}d(t,T_{n_i})\ls \\
&\ls 2(\sum^{\infty}_{i=0}2^{-i})d(t,T_m)\ls 4d(t,T_m)
\end{align*}
and $\sum^{\infty}_{i=0}\frac{1}{N_{n_i}^3}V_{n_i}\ls Z$. Thus we finally get
$$
|X(t)-X(\pi_m(t))|\ls 15K^2\sigma_{m}(t)+4K^2d(t,T_m)\psi^{-1}(Z).
$$
It completes the proof with $A=15K^2$ and $B=4K^2$.
\end{dwd}
\begin{rema}
Note that the above proof works for any $\psi$ that is not necessarily convex, 
however then the constants $A$, $B$ may depend on $\psi$ since there may be no longer true that $\psi(2^i)\gs 2^i$. 
In particular for $\psi(x)=\varphi_p(x)=2^{x^p}-1$, $0<p\ls 1$ we have
$K=2^{\frac{1}{p}-1}$ and  $\E Z\ls (2^{\frac{1}{p}}-1)^{-1}$. Consequently
$Z$ has to be multiplied by $2^{\frac{1}{p}}-1$ which affects $A,B$. 
\end{rema}
\begin{coro}
Under assumptions of Proposition \ref{pro1}, there exist universal constants $A,B$ such that
$$
\E \sup_{m\gs 0}\sup_{t\in T}\psi(\frac{(|X(t)-X(\pi_m(t))|-A\sigma_m(t))_{+}}{Bd(t,\pi_m(t))})\ls 1. 
$$ 
\end{coro}
\smallskip

\noindent Now we turn to analyze the modulus of continuity of $X(t)$, $t\in T$. Consider any $s,t\in T$.
There exists the smallest $k\gs 0$ such that both $d(s,T_k)$ and $d(t,T_k)$ are less than $d(s,t)$, i.e.
$$
k(s,t)=\max\{k\gs 0:\;d(s,T_k)+d(t,T_k)\gs d(s,t)\}.
$$
Note that $k(s,t)$ is well defined since at least $d(s,T_0)+d(t,T_0)\gs d(s,t)$.
Obviously $k(s,t)+1$ is the required level where in the chaining construction it is better to jump from 
the approximation of $t$ to the approximation of $s$. We explain proving that $k(s,t)$ the value $k\gs 0$
for which the function
$$
f(k)=\sigma_{k+1}(s)+\sigma_{k+1}(t)+(2^{k+1}-1)d(s,t)
$$
is the smallest possible. Indeed $f(k)> f(k-1)$ implies that
$$
2^k(d(s,T_k)+d(t,T_k))> 2^k d(s,t)
$$ 
and hence $k> k(s,t)$. On the other hand for all $k\ls k(s,t)$, we have that
$f(k)\ls f(k-1)$. Therefore $k(s,t)$ is the argument minimum of $f$.
Consequently in the view of the proof of Proposition \ref{pro1} the deterministic part
that bounds $|X(t)-X(s)|$ should be up to a constant bounded by 
$$
\bar{\tau}(s,t)=f(k(s,t))=\sigma_{k(s,t)+1}(s)+\sigma_{k(s,t)+1}(t)+2^{k(s,t)+1}d(s,t).
$$
We are ready to define a simple distance $\tau$ on $T$ such that $\tau(s,t)$ is comparable with this quantity. 
\smallskip

\noindent
Define the following numbers 
$$
\sigma(t,a)=\sum^{\infty}_{n=0}2^{n}\min\{d(t,T_n),a\},\;\;a> 0,\;\;t\in T.
$$
Let $\tau$ be a new distance on $T$ given by
$$
\tau(s,t)=\max\{\sigma(t,d(s,t)),\sigma(s,d(s,t))\},\;\;\mbox{for all}\;s,t\in T.
$$
We show that $\tau(s,t)$ is indeed a distance on $T$ i.e. it satisfies the triangle inequality.
It suffices to check that for $s,t,u\in T$ there holds
$$
\tau(s,u)\ls \tau(s,t)+\tau(t,u).
$$
First note that
$$
\tau(s,u)\ls \max\{\sigma(s,d(s,t)+d(t,u)),\sigma(u,d(s,t)+d(t,u))\}.
$$
and then 
\begin{align*}
&\sigma(s,d(s,t)+d(t,u))\ls \sum^{\infty}_{n=0}2^{n}d(s,T_n)1_{d(s,T_n)\ls d(s,t)}+\\
&+\sum^{\infty}_{n=0}2^{n}\min\{d(s,t)+d(t,T_n),d(s,t)+d(t,u)\}1_{d(s,T_n)>d(s,t)}\ls\\
&\ls \sum^{\infty}_{n=0}2^{n}\min\{d(s,t),d(s,T_n)\}+\sum^{\infty}_{n=0}2^{n}\min\{d(t,u),d(t,T_n)\}\ls \tau(s,t)+\tau(t,u).
\end{align*}
In the same way we get
$$
\sigma(u,d(s,t)+d(t,u))\ls \tau(s,t)+\tau(t,u).
$$
The distance $\tau(s,t)$ is comparable with $\bar{\tau}(s,t)$, namely
\be\label{megatron}
2^{-1}\bar{\tau}(s,t)\ls\tau(s,t)\ls \bar{\tau}(s,t).
\ee
Indeed since for $k=k(s,t)$ we have $d(s,T_l)+d(t,T_l)\gs d(s,t)$ for $l\ls k$ it implies that
\begin{align*}
& \tau(s,t)\gs \max\{\sigma_{k+1}(s)+\sum^k_{l=0}2^{k}(d(s,T_k)\wedge d(s,t)),
\sigma_{k+1}(t)+\sum^k_{l=0}2^l(d(t,T_k)\wedge d(s,t))\}\gs \\
&\gs  2^{-1}[\sigma_{k+1}(s)+\sigma_{k+1}(t)+2(2^{k+1}-1)d(s,t)]=2^{-1}\bar{\tau}(s,t).
\end{align*}
On the other hand
\begin{align*}
& \tau(s,t)\ls \max\{\sigma_{k+1}(s),\sigma_{k+1}(t)\}+\sum^k_{l=0} 2^l d(s,t)]\ls \\
& \ls \sigma_{k+1}(s)+\sigma_{k+1}(t)+2^{k+1}d(s,t)=\bar{\tau}(s,t).
\end{align*}
We recall that
$$
Z=\sum^{\infty}_{n=1}\frac{1}{N_n^3}\sum_{t\in T_{n}}\psi(\frac{|X(t)-X(\pi_{n-1}(t))|}{d(t,\pi_{n-1}(t))}).
$$
Let us state the main result of this section that $\tau(s,t)$ is the 
suitable modulus of continuity.
\begin{prop}\label{pro2}
For each admissible net $(T_n)^{\infty}_{n=0}$, and separable $X(t)$, $t\in T$ that satisfies
(\ref{erm1}) the following inequality holds. For all $s,t\in T$
$$
|X(s)-X(t)|\ls \bar{A}\tau(s,t)+\bar{B}d(s,t)\psi^{-1}(Z),
$$
where $\bar{A}=30K^2,\bar{B}=10K^2$, $Z\gs 0$ and $\E Z\ls 1$.
\end{prop}
\begin{dwd}
Let $k=k(s,t)$ for $s,t\in T$. Observe that
\begin{eqnarray}
&& |X(t)-X(s)|\ls |X(t)-X(\pi_{k+1}(t))|+\nonumber \\
\label{zero} && +|X(s)-X(\pi_{k+1}(t))|+|X(\pi_{k+1}(t))-X(\pi_{k+1}(s))|.
\end{eqnarray}
By Proposition \ref{pro1} we have that for $x\in\{s,t\}$
$$
|X(x)-X(\pi_{k+1}(x))|\ls A\sigma_{k+1}(x)+Bd(x,\pi_{k+1}(x))\psi^{-1}(Z), 
$$
where $A=15K^2$ and $B=4K^2$. Since 
$$
d(s,\pi_{k+1}(s))+d(t,\pi_{k+1}(t))\ls d(s,t)
$$
it implies that
\begin{eqnarray}\label{jeden}
&& |X(s)-X(t)-X(\pi_{k+1}(s))+X(\pi_k(t))|\ls \nonumber\\
&& \ls  A(\sigma_{k+1}(s)+\sigma_{k+1}(t))+Bd(s,t)\psi^{-1}(Z).
\end{eqnarray}
On the other hand using our main argument
\begin{align*}
& |X(\pi_{k+1}(t))-X(\pi_{k+1}(s))|\ls \\
&\ls d(\pi_{k+1}(t),\pi_{k+1}(s))(3K^2 2^{k+1}+K^2\psi^{-1}(\frac{1}{N_{k+1}^3}V_{k+1})),
\end{align*}
where
$$
V_{k+1}=\sum_{u,v\in T_{k+1}}\varphi(\frac{|X(u)-X(v)|}{d(u,v)}).
$$
Note that by the the definition of $k(s,t)$
$$
d(\pi_{k+1}(t),\pi_{k+1}(s))\ls d(t,T_{k+1})+d(s,t)+d(s,T_{k+1})\ls 2d(s,t).
$$
Therefore using that $\sum^{\infty}_{k=0}\frac{1}{N_{k+1}^3}V_{k+1}\ls Z$
$$
|X(\pi_{k+1}(s))-X(\pi_{k+1}(t))|\ls 6K^2 2^{k+1}d(s,t)+2K^2d(s,t)\varphi^{-1}(Z).
$$
and hence 
\be\label{dwa}
|X(\pi_{k+1}(s))-X(\pi_{k+1}(t))|\ls 15K^2(2^{k+1}-1)d(s,t)+2K^2d(s,t)\varphi^{-1}(Z).
\ee
Now we have to sum up bounds (\ref{zero}), (\ref{jeden}), (\ref{dwa})
and apply (\ref{megatron}). The proof is completed with $\bar{A}=30K^2$ and $\bar{B}=10K^2$.
\end{dwd}
\begin{coro}
The following inequality holds
$$
\E \sup_{s,t\in T}\psi((\frac{|X(s)-X(t)|-A\tau(s,t))_{+}}{Bd(s,t)})\ls 1,
$$
where $A,B$ are constant for Proposition \ref{pro2}.
\end{coro}
The meaning of the result is that for exponentially concentrated random variables there is always some form Law of Iterated Logarithm. For example on the interval $[0,1]\subset\R$
and any fractional Brownian motion $X(t)$, $t\in T$ with the Hurst exponent $H\in (0,1)$,
we have $d(s,t)=|s-t|^{H}$ and with usual $T_n=\{kN_n^{-1}:\; k\in \{1,...,N_n\}\}\}$
we get $\tau(s,t)\sim |s-t|^H\sqrt{\log_2(1+|s-t|^{-H})}$.
\begin{coro}
If $X(t)$, $t\in [0,1]$ is a fractional Brownian motion with the Hurst coefficient $H\in (0,2)$ then
there exists universal $A,B$ such that
$$
\E\sup_{s,t\in [0,1]}\varphi_2(\frac{(|X(s)-X(t)|-A |s-t|^H\sqrt{\log_2(1+|s-t|^{-H})})_{+}}{B|s-t|^{H}})\ls 1.
$$
\end{coro}

\section{The concentration inequality for two distances}

As the basic example of the Bernstein inequality shows one distance may not suffice
to fully describe the concentration property. The classical situation
when this happens concerns independent, symmetric, identically distributed $X_i$, $1\ls i\ls n$ 
of log concave tails, i.e. when $u\ra -\log \P(|X_i|>u)$ is a convex function for $u\gs 0$.
We can consider canonical type process using $T\subset \R^n$ and  $X(t)=\sum^n_{i=1}t_i X_i$.
Therefore we generalize slightly the idea described 
in the previous sections towards the case where two distances can applied.
\smallskip

\noindent Assume that on set $T$ there are two distances $d_1,d_2$ that imply the same 
same topology on $T$ (i.e. we require that the convergence in $d_1$ is equivalent to the convergence in $d_2$).
The main assumption on increments, is that there exists Young functions $\psi_1$ and $\psi_2$
such that 
\be\label{erm4}
\E\min\{\psi_1(\frac{|X(s)-X(t)|}{d_1(s,t)}),\psi_2(\frac{|X(s)-X(t)|}{d_2(s,t)})\}\ls 1 \;\;\mbox{for} \;s,t\in T.
\ee
The same remark as for one distance control is valid. We can extend the approach on $\psi_1$, $\psi_2$
that are increasing to infinity, continuous and $\psi_i(0)=0$, $\psi_i(1)=1$ for $i=1,2$.
Moreover we assume the condition (\ref{erm3}) for $\psi_1$ and $\psi_2$, yet it is not enough for the our analysis. We need a polynomial comparability of $\psi^{-1}_1$ and $\psi^{-1}_2$. We state
the condition in the general form but to avoid technical complications we assume that there exists a single $\psi$ that satisfies (\ref{erm3})  such that $\psi_1(x)=\psi(x^{p_1})$, and $\psi_2(x)=\psi(x^{p_2})$ for $p_1,p_2>0$. The simplest case of the setting is when $\psi_1=\varphi_{p_1}$ and $\psi_2=\varphi_{p_2}$ for some $p_1,p_2\gs 1$
since then $\psi_1$ and $\psi_2$ are convex.
\smallskip

\noindent
For simplicity let us assume that $T$ is compact in the topology defined by $d_1$ and $d_2$.
Moreover we require that $\psi$ such that $\psi_1(x)=\psi(x^{p_1})$ and $\psi_2(x)=\psi(x^{p_2})$ are convex.
\smallskip

\noindent 
The main point is that we define the common approximation net for both two distances. 
Let $N_n=\psi(2^{n})$, for $n\gs 1$ and $N_0=1$. Again we assume that $(T_n)^{\infty}_{n=0}$ is admissible, i.e. $|T_n|\ls N_{n}$, $T_{n}\subset T_{n+1}$ and $\bigcup^{\infty}_{n=0}T_n$
is dense in $T$. Consequently $\psi^{-1}_1(N_n)=2^{p_1 n}$ and $\psi^{-1}_2(N_n)=2^{p_2 n}$.
We extend our definition of $\sigma_{m}(t)$, i.e. we define
$$
\sigma^{j}_{m}(t)=\sum^{\infty}_{n=m}2^{p_j n}d_j(t,T_n),\;\;j\in{1,2}.
$$
W.l.o.g. we may assume that $\sigma^j(t)<\infty$ for $j=1,2$.
Moreover let 
$$
Z=\sum^{\infty}_{n=1}\frac{1}{N_n^3}\sum_{t\in T_{n}}\min\{\psi_1(\frac{|X(t)-X(\pi_{n-1}(t))|}{d_1(t,\pi_{n-1}(t))}),
\psi_2(\frac{|X(t)-X(\pi_{n-1}(t))|}{d_2(t,\pi_{n-1}(t))})\}.
$$
Fix $t\in T$, we aim to extend Proposition \ref{pro1}. Obviously again for a given $t\in T$ we can use the sequence 
$n_i=n_i(t)$ such that $n_0=m$ and it increases to $\infty$, but the definition is  more complicated.
Let $n_0=m$ and 
\begin{align*}
& n_i=\inf\{n>n_{i-1}:\;2^{p_1 n}d_1(t,T_{n_{i-1}})+2^{p_2 n}d_2(t,T_{n_{i-1}})>\\
& >2(2^{p_1 n}d_1(t,T_n)+2^{p_2 n}d_2(t,T_n))\}.
\end{align*}
In particular it means that either $d_1(t,T_{n_{i-1}})>2 d_1(t,T_{n_i})$
or $d_2(t,T_{n_{i-1}})>2 d_2(t,T_{n_i(t)})$. However we cannot claim that the property holds for both two distances at once.
It should be noticed also why $n_i$ necessarily exists and this is due to the assumption $\sigma_j(t)<\infty$ for $j\in \{1,2\}$.
\smallskip

\noindent
To formulate the best possible result we define for each $t\in T$
$$
\bar{d}_1(t,\pi_m(t))=\sum^{\infty}_{i=0}d_1(t,T_{n_i}),\;\;\bar{d}_2(t,\pi_m(t))=\sum^{\infty}_{i=0}d_2(t,T_{n_i}).
$$ 
Obviously in if $T_n$ well approximates $T$ in both two distances then $\bar{d}_j(t,\pi_m(t))$ is comparable with $d_j(t,\pi_m(t))$ but in full generality we cannot claim the statement.
On the other hand there are estimates on $\bar{d}_j(t,\pi_m(t))$ in terms of $\sigma^j_m(t)$, i.e.
$$
\bar{d}_j(t,\pi_m(t))\ls 2^{-p_j m}\sigma^j_m(t),\;\;\mbox{for}\;j\in \{1,2\}.
$$ 
Usually the above bounds are sufficiently strong for concentration reason. 
\begin{prop}\label{pro3}
Let $X(t)$, $t\in T$ satisfy (\ref{erm4}) with $\psi_1=\psi(x^{p_1})$ and $\psi_2=\psi(x^{p_2})$ which verify
\ref{erm3} with constants $K_1,K_2$ respectively. 
For all $t\in T$ and $m\gs 0$ the following inequality holds
\begin{align*}
& |X(t)-X(\pi_m(t))|\ls A(\sigma^1_{m}(t)+\sigma^2_{m}(t))+ \\
& +B(\bar{d}_1(t,\pi_m(t))\psi^{-1}_1(Z)+\bar{d}_2(t,\pi_m(t))\psi^{-1}_2(Z)),
\end{align*}
where $A=3(1+2^{1+p})K^2,B=2K^2$ and $K=\max\{K_1,K_2\}$, $p=\max\{p_1,p_2\}$.
\end{prop}
\begin{dwd}
We use the sequence $n_i=n_i(t)$, $i\gs 0$ to get
\be\label{a22}
|X(t)-X(\pi_m(t))|\ls \sum^{\infty}_{i=0}|X(\pi_{n_i}(t))-X(\pi_{n_{i+1}}(t))|.
\ee
Consequently due to (\ref{erm3})
\begin{align*}
& |X(\pi_{n_i}(t))-X(\pi_{n_{i-1}}(t))|\ls\\
&\ls \sum_{j\in\{1,2\}} d_j(\pi_{n_{i}}(t),\pi_{n_{i-1}}(t))\psi^{-1}_j(\frac{N_{n_{i}}^3}{N_{n_{i}}^3}V_{n_i})\ls \\
&\ls  \sum_{j\in \{1,2\}}d_j(\pi_{n_{i}}(t),\pi_{n_{i-1}}(t))[3K^2_j 2^{p_j n_{i}}+K^2_j\psi^{-1}_j(\frac{1}{N_{n_{i}}^3}\sum_{u,v\in T_{n_i}}V_{n_i})],
\end{align*}
where
$$
V_n=\sum_{u,v\in T_{n}}\min\{\psi_1(\frac{|X(u)-X(u)|}{d_1(u,v)}),\psi_2(\frac{|X(u)-X(v)|}{d_2(u,v)})\}\;\;\mbox{for}\;n\gs 0.
$$
Observe that for $j\in\{1,2\}$
$$
d_j(\pi_{n_{i}}(t),\pi_{n_{i-1}}(t))\ls d_j(t,T_{n_i})+d_j(t,T_{n_{i-1}})
$$
and therefore by the construction of $n_i$
\begin{align*}
& \sum_{j\in \{1,2\}}3K^2_j  d_j(\pi_{n_{i}}(t),\pi_{n_{i-1}}(t)) 2^{p_j n_{i}}\ls 
3K^2\sum_{j\in \{1,2\}}(d_j(t,T_{n_i})+d(t,T_{n_{i-1}}))2^{p_j n_i}\ls \\ 
& \ls 3K^2 \sum_{j\in \{1,2\}}(d(t,T_{n_i})+2d(t,T_{n_i-1}))2^{p_jn_i}.
\end{align*}
Consequently
$$
\sum^{\infty}_{i=0 }\sum_{j\in \{1,2\}}3K^2_j d_j(\pi_{n_{i}}(t),\pi_{n_{i-1}}(t)) 2^{p_j n_{i}}\ls 3(1+2^{1+p})K^2(\sigma^1_m(t)+\sigma^2_m(t)). 
$$
On the other hand
$$
\sum^{\infty}_{i=0}K^2_jd_j(\pi_{n_{i}}(t),\pi_{n_{i-1}}(t))\ls 2K^2_j\bar{d}_j(t,\pi_m(t)).
$$
and $\sum^{\infty}_{i=0}\frac{1}{N_{n_i}^3}V_{n_i}\ls Z$, hence
\begin{align*}
& |X(t)-X(\pi_m(t))|\ls 3(1+2^{1+p})K^2(\sigma^1_m(t)+\sigma_m^2(t))+\\
&+2K^2(\bar{d}_1(t,\pi_m(t))\psi_1^{-1}(Z)+\bar{d}_2(t,\pi_m(t))\psi_2^{-1}(Z)).
\end{align*}
It completes the proof with $A=3(1+ 2^{1+p})K^2$ and $B=2K^2$, where $K=\max\{K_1,K_2\}$ and $p=\max\{p_1,p_2\}$.
\end{dwd}
As for the modulus of continuity for given $s,t\in T$ one has to be more careful. 
For the one distance control it was clear for which $k\gs 0$ it is worth to jump from the approximation of  $t$ to the approximation of $s$. For the two distance control we have two possible solutions $k^1(s,t)$ and $k^2(s,t)$
such that for $j\in\{1,2\}$
$$
k^j(s,t)=\max\{k\gs 0:\; \sum_{x\in\{s,t\}}d_j(x,T_k)\gs 
d_j(s,t)\}.
$$
Using that $p_1,p_2\gs 1$ can prove in the same way as for the one distance control that
$k^j(s,t)$, $j\in\{1,2\}$ which is the $k\gs 0$ for which the function
$$
\sigma^j_{k+1}(s)+\sigma^j_{k+1}(t)+\frac{2^{p_j(k+1)}-1}{2^{p_j}-1}d_j(s,t)
$$
is the smallest possible. We use the same idea to define the general $k(s,t)$ in the case we study.
\smallskip

\noindent
First define $k(s,t)$ as
$$
k(s,t)=\max\{k\gs 0:\; \sum_{j\in \{1,2\}}\sum_{x\in \{s,t\}}
2^{p_j k}d_j(x,T_k)\gs \sum_{j\in\{1,2\}}2^{p_j k}d_j(s,t)\}.
$$
As above we explain that for such $k=k(s,t)$ the function
$$
f(k)=\sum_{x\in\{s,t\}}\sum_{j\in \{1,2\}}\sigma^j_{k+1}(x)+
\sum_{j\in \{1,2\}} \frac{2^{p_j(k+1)}-1}{2^{p_j}-1}d_j(s,t),\;\;k\gs 0 
$$
is the smallest possible. Indeed observe that $f(k)> f(k-1)$ implies that
$$
\sum_{j\in \{1,2\}}2^{p_jk} d_j(s,t)\gs \sum_{j\in \{1,2\}}\sum_{x\in \{s,t\}}2^{p_j k}d_j(x,T_k), 
$$
which happens for all $k> k(s,t)$. On the other hand for all $k\ls k(s,t)$
we have that $f(k)\ls f(k-1)$. 
$$
\sum_{j\in \{1,2\}}2^{p_jk}d_j(s,t)\gs \sum_{j\in \{1,2\}}
\sum_{x\in \{s,t\}}2^{p_jk}d_j(x,T_k), 
$$
By the property of $k^j(s,t)$, $j\in \{1,2\}$ we have mentioned above
we get
\be\label{miau1}
\min\{k^1(s,t),k^2(s,t)\} \ls k(s,t)\ls \max\{k^1(s,t),k^2(s,t)\}.
\ee
We use $k(s,t)$ to define the deterministic bound on $|X(t)-X(s)|$. Let
$$
\bar{\tau}(s,t)=\sum_{x\in \{s,t\}}\sum_{j\in \{1,2\}}\sigma^j_{k(s,t)+1}(x)+\sum_{j\in \{1,2\}} \frac{2^{p_j(k(s,t)+1)}-1}{2^{p_j}-1}d_j(s,t).
$$
THe obvious result is
\be\label{owoc}
\tau_1(s,t)+\tau_2(s,t)\ls \bar{\tau}(s,t),
\ee
where
$$
\tau_j(s,t)=\max\{\sigma_j(t,d(s,t)),\sigma_j(s,d(s,t))\},\;\;\mbox{for}\;s,t\in T, \;j\in \{1,2\}
$$ 
and
$$
\sigma_j(t,a)=\sum^{\infty}_{n=0}2^{p_j n}\min\{d_j(t,T_n),a\},\;\;a> 0,\;\;t\in T.
$$
Indeed for $k=k(s,t)$ and $j\in \{1,2\}$
\begin{align*}
& \tau_j(s,t)\ls \max_{x\in \{s,t\}} \sigma^j_{k+1}(x)+\sum^k_{l=0} 2^{p_j l}d_j(s,t)\ls \\
&\ls \sigma^j_{k+1}(s)+\sigma^j_{k+1}(t)+\frac{2^{p_j(k+1)}-1}{2^{p_j}-1}d_j(s,t)
\end{align*}
and hence
$$
\tau_1(s,t)+\tau_2(s,t)\ls \bar{\tau}(s,t).
$$
Unfortunately $\tau_1(s,t)+\tau_2(s,t)$ may be not comparable with $\bar{\tau}(s,t)$.
Nevertheless we can still provide an upper bound, let $k_1(s,t)=k_1$, $k_2(s,t)=k_2$ the following holds
\begin{eqnarray}\label{gus}
&& \bar{\tau}(s,t)\ls 2(\tau_1(s,t)+\tau_2(s,t))+\frac{2^{p_2}}{2^{p_2}-1}(2^{k_1p_2}-2^{k_2 p_2})d_2(s,t)1_{k_1\gs k_2}+\nonumber\\
&&+\frac{2^{p_1}}{2^{p_1}-1}(2^{k_2 p_1}-2^{k_1 p_1})d_1(s,t)1_{k_2\gs k_1}.
\end{eqnarray}
Assume for simplicity that $k_1\ls k_2$. By the definition of $k(s,t)$
\begin{align*}
& \bar{\tau}(s,t)\ls \sigma^1_{k_2+1}(s)+\sigma^1_{k_2+1}(t)+\frac{2^{(k_2+1)p_1}-1}{2^{p_1}-1}d_1(s,t)+\\
&+\sigma^2_{k_2+1}(s)+\sigma^2_{k_2+1}(t)+\frac{2^{(k_2+1)p_2}-1}{2^{p_2}-1}d_2(s,t).
\end{align*}
Now observe that
$$
\sigma^2_{k_2+1}(s)+\sigma^2_{k_2+1}(t)+\frac{2^{(k_2+1)p_2}-1}{2^{p_2}-1}d_2(s,t)\ls 2\tau_2(s,t)
$$
and using that $k_1\ls k_2$ 
\begin{align*}
& \sigma^1_{k_2+1}(s)+\sigma^1_{k_2+1}(t)+\frac{2^{(k_2+1)p_1}-1}{2^{p_1}-1}d_1(s,t)\ls \sigma^1_{k_1+1}(s)+\sigma^1_{k_1+1}(t)+\\
&+\frac{2^{(k_1+1)p_1}-1}{2^{p_1}-1}d_1(s,t)+\frac{2^{(k_2+1)p_1}-2^{(k_1+1)p_1}}{2^{p_1}-1}d_1(s,t)\ls\\
&\ls 2\tau_1(s,t)+\frac{2^{p_2}}{2^{p_2}-1}(2^{k_2p_1}-2^{k_1p_1})d_1(s,t).
\end{align*}
It implies (\ref{gus}). To understand this inequality better we have to observe that
$$
2^{k_2p_1}\ls (2^{k_2p_2})^{\frac{p_1}{p_2}}\ls \left(\frac{\tau_2(s,t)}{d_2(s,t)}\right)^{\frac{p_1}{p_2}}
$$
and hence
$$
\frac{2^{p_1}}{2^{p_1}-1}(2^{k_2p_1}-2^{k_1p_1})d_1(s,t)\ls \frac{2^{p_2}}{2^{p_1}-1}\left(\frac{\tau_2(s,t)}{d_2(s,t)}\right)^{\frac{p_1}{p_2}}.
$$
Consequently
\begin{align*}
&\bar{\tau}(s,t)\ls \tau_1(s,t)+\tau_2(s,t)+\\
&+\max\{\frac{2^{p_2}}{2^{p_2}-1}\left(\frac{\tau_2(s,t)}{d_2(s,t)}\right)^{\frac{p_1}{p_2}}d_1(s,t),\frac{2^{p_2}}{2^{p_2}-1}\left(\frac{\tau_1(s,t)}{d_1(s,t)}\right)^{\frac{p_2}{p_1}}d_2(s,t) \}
\end{align*}
We need also a slight generalization of the distances $d_1$ and $d_2$, namely
$$
\bar{d}_j(s,t)=\sum^{\infty}_{i=0}d_j(t,T_{n_i(t)})+d_j(s,t)+\sum^{\infty}_{i=0}d_j(s,T_{n_i(s)}),
$$
and $n_i(t), n_i(s)$ are sequences $n_i$ for $t,s$ (recall that $n_i$ depend on points in $T$) that starts from $k(s,t)+1$. Again observe that
\be\label{rar2}
\bar{d}_1(s,t)+\bar{d}_2(s,t)\ls \max
\{2^{-p_1k(s,t)},2^{-p_2k(s,t)}\} \bar{\tau}(s,t).
\ee
We may state the main result of these section.
\begin{prop}\label{pro4}
For each admissible net $(T_n)^{\infty}_{n=0}$, and separable $X(t)$, $t\in T$ that satisfies
(\ref{erm1}) the following inequality holds. For all $s,t\in T$
$$
|X(s)-X(t)|\ls \bar{A}(\bar{\tau}(s,t))+\bar{B}[\bar{d}_1(s,t)\psi^{-1}_1(Z)+\bar{d}_2(s,t)\psi^{-1}_2(Z)],
$$
where $\bar{A}=3(1+2^{1+p})K^2,\bar{B}=5K^2$, $Z\gs 0$ and $\E Z\ls 1$.
\end{prop}
\begin{dwd}
Let $k=k(s,t)$ for $s,t\in T$. Observe that
\ba\label{zero2}
&& |X(t)-X(s)|\ls |X(t)-X(\pi_{k+1}(t))|+\nonumber \\
&&+|X(s)-X(\pi_{k+1}(t))|+|X(\pi_{k+1}(t))-X(\pi_{k+1}(s))|.
\ea
By Proposition \ref{pro3}  we have for $x\in\{s,t\}$
\begin{eqnarray}
&& |X(x)-X(\pi_{k+1}(x))|\ls A[\sigma^1_{k+1}(x)+\sigma^2_{k+1}(x)]+\nonumber\\
\label{jeden2}&&+B[\bar{d}_1(x,\pi_k(x))\psi^{-1}_1(Z)+\bar{d}_2(x,\pi_k(x))\psi^{-1}_2(Z)], 
\end{eqnarray}
where $A=3(1+2^{1+p})K^2$ and $B=2K^2$.
We use (\ref{erm3}) for $\psi_1$ and $\psi_2$ with constants $K_1,K_2$
\begin{align*}
& |X(\pi_{k+1}(t))-X(\pi_{k+1}(s))|\ls \\
&\ls \sum_{j\in \{1,2\}}d_j(\pi_{k+1}(t),\pi_{k+1}(s))(3K^2_j2^{p_j(k+1)}+K^2_j\psi_j^{-1}(\frac{1}{N_{k+1}^3}V_{k+1})).
\end{align*}
where
$$
V_{k+1}=\sum_{u,v\in T_{k+1}}
\min\{\psi_1(\frac{|X(u)-X(v)|}{d_1(u,v)}),\psi_2(\frac{|X(u)-x(v)|}{d_2(u,v)})\}.
$$
Note that 
$$
d_j(\pi_{k+1}(t),\pi_{k+1}(s))\ls d_j(t,T_{k+1})+d_j(s,t)+d_j(s,T_{k+1}).
$$
Therefore by the construction of $k$
$$
\sum_{j\in \{1,2\}}K^2_j 2^{p_j (k+1)}d_j(\pi_{k+1}(s),\pi_{k+1}(t))\ls  2K^2\sum_{j\in \{1,2\}}2^{p_j (k+1)}d_j(s,t).
$$
Consequently
\begin{align*}
&|X(\pi_{k+1}(s))-X(\pi_{k+1}(t))|\ls 6K^2 \sum_{j\in \{1,2\}}2^{p_j (k+1)}d_j(s,t) +\\
&+K^2[\bar{d}_1(s,t)\psi^{-1}_1(Z)+\bar{d}_2(s,t)\psi^{-1}_2(Z)].
\end{align*}
Therefore
\begin{eqnarray}\label{dwa2}
&& |X(\pi_{k+1}(s))-X(\pi_{k+1}(t))|\ls 6K^2 2^p \sum_{j\in \{1,2\}}\frac{2^{p_j(k+1)}-1}{2^{p_j}-1}d_j(s,t)+\nonumber\\
&&+K^2[\bar{d}_1(s,t)\psi^{-1}_1(Z)+\bar{d}_2(s,t)\psi^{-1}_2(Z)].
\end{eqnarray}
Now we have to sum up bounds (\ref{zero2}), (\ref{jeden2}), (\ref{dwa2})
$$
|X(s)-X(t)|\ls 2A\bar{tau}(s,t)+5K^2[\bar{d}_1(s,t)\psi^{-1}_1(Z)+\bar{d}_2(s,t)\psi^{-1}_2(Z)].
$$
The proof is completed with $\bar{A}=6(1+2^{1+p})K^2$ and $B=5K^2$.
\end{dwd}
The consequence of the theorem is the Bernstein type inequality for chaining. 
Using (\ref{rar2}) and Proposition \ref{pro3}.
\begin{coro}
The following inequality holds
$$
\E\min\{\psi_1(\frac{|X(s)-X(t)|-\bar{A}\bar{\tau}(s,t)}{\bar{B}\bar{d}_1(s,t)}),
\psi_2(\frac{|X(s)-X(t)|-\bar{A}\bar{\tau}(s,t)}{\bar{B}\bar{d}_2(s,t)})\}\ls 1.
$$
Moreover 
$$
\bar{d}_j(s,t)\ls 2^{-p_j(k(s,t)+1)}\tau_j(s,t),\;\;\mbox{for}\;j\in \{1,2\},
$$
where $k(s,t)$ is defined by (\ref{rar2}).
\end{coro}
Obviously we natural application of the idea should be to the theory of empirical processes.
We discuss the question for a certain question in the following section.

\section{Square estimation for one distance}

The following problem was studied in \cite{MPT}.
Let $X_1,X_2,...,X_N$ be independent random variables with values in a measurable space $(\ccX,\ccB)$.
Let $\ccF$ be a family of real measurable functions on $(\ccX,\ccB)$. To explain the problem of square estimation
we start from the analysis of one distance control and then we turn to consider more complicated case of two distance control.
\smallskip

\noindent
Let us assume that there exists distance $d$ on $\ccF$ such that for any given $f,g\in \ccF$
\be\label{eq:5.1}
\P(|(f-g)(X_i)|\gs u)\ls 2\exp(-\frac{u^2}{d(f,g)^2}),
\ee
for all $1\ls i\ls N$, $u>0$. Moreover we assume that
\be\label{eq:5.2}
\P(|f(X_i)|\gs u)\ls 2\exp(-\frac{u^2}{d(f,0)^2})
\ee
Recall that it is equivalent to $\|(f-g)(X_i)\|_{\varphi_2}\ls Cd(f,g)$ and $\|f(X_i)\|_{\varphi_2}\ls Cd(f,0)$ 
for all $i\in \{1,2,...,N\}$ and a universal constant $C>0$. Obviously in particular
\be\label{eq:5.3}
\E (f-g)^2(X_i)\ls 2C^2d(f,g)^2\;\;\mbox{and}\;\;\frac{1}{N}\E \sum^N_{i=1}(f-g)^2(X_i)\ls 2C^2d(f,g)^2
\ee
and similarly
\be\label{eq:5.4}
\E f^2(X_i)\ls 2C^2d(f,0)^2\;\;\mbox{and}\;\;\frac{1}{N}\E \sum^N_{i=1}f^2(X_i)\ls 2C^2d(f,0)^2
\ee
We aim to provide a concentration inequality for
$$
S_{N}(f)=\frac{1}{N}\sum^N_{i=1}(f^2-\E f^2)(X_i).
$$
To state the result we have to recall Bernstein type inequalities.
For all $f,g\in \ccF$ the following holds
\be\label{eq:5.45}
\P(|S_N(f-g)|\gs u)\ls 2\exp(-N\min(\frac{u^2}{4
d(f,g)^4},\frac{u}{4d(f,g)^2}))
\ee
Moreover for all $f\in \ccF$
$$
\P(|S_N(f)|\gs u)\ls 2\exp(-N\min(\frac{u^2}{4d(f,0)^4}),\frac{u}{4d(f,0)^2}).
$$
The above inequalities can be rewritten using the following function
$$
\varphi(x)=2^{\min\{\sqrt{N}x,x^2\}}-1,\;\;\varphi^{-1}(x)=\max\{\frac{1}{\sqrt{N}}
\log_2(1+x),\sqrt{\log_2(1+x)}\}.
$$
The function $\varphi$ is not convex but is comparable to a convex function
say $\psi(x)=2^{\min\{2\sqrt{N}x-N,x^2\}}$ (note that $\psi'(\sqrt{N})=2\log 2\sqrt{N}2^N$) . Clearly
$$
\psi(x/2)\ls \varphi(x)\ls \psi(x),\;\;\mbox{for all}\;\;x>0.
$$
As usual  $\varphi(0)=0$ and $\varphi(1)=1$ and $\varphi$ and satisfies (\ref{erm3}) with $K=1$. Moreover $\varphi(D^{-1}x)\ls D^{-1}\varphi(x)$ for $D\gs 1$ and $x\gs 0$. 
The meaning of (\ref{eq:5.1}) in terms of $\varphi$ is that $f,g\in \ccF$
\be\label{burak3}
\|S_N(f-g)\|_{\varphi}\ls K^2 N^{-\frac{1}{2}}d(f,g)^2.
\ee
In the same way the meaning of (\ref{eq:5.2}) is such that for each $f\in \ccF$
\be\label{burak4}
\|S_N(f)\|_{\varphi}\ls K^2 N^{-\frac{1}{2}}d(f,0)^2.
\ee
On the other hand
$$
\| f^2-g^2\|_{\varphi_1}\ls \|f-g\|_{\varphi_2}\|f+g\|_{\varphi_2}\ls C^2(d(f,0)+d(g,0))d(f,g).
$$
Consequently for any $f,g\in \ccF$
\be\label{eq:5.5}
\| S_N(f)-S_N(g)\|_{\varphi}\ls K^2 N^{-\frac{1}{2}}(d(f,0)+d(g,0))d(f,g).
\ee
We are ready to state the main result of the section.
\begin{theo}\label{thm1}
Suppose that $\alpha=\sup_{t\in \ccF}d(f,0)<\infty$ for all $f\in \ccF$  then
$$
\sup_{f\in \ccF} \frac{1}{N}\sum^N_{i=1}(f^2-\E f^2)(X_i)\ls A(N^{-1}\gamma^2_2(\ccF)+N^{-\frac{1}{2}}\alpha\gamma_2(\ccF))+BN^{-\frac{1}{2}}\alpha^2\varphi^{-1}(Z),
$$
where $Z\gs 0$ is such that $\E Z\ls 1$ and $A,B$ are universal constants. In terms of concentration it means that for all $u\gs 0$
\begin{align*}
&\P(\sup_{f\in \ccF}\frac{1}{N}(\sum^N_{i=1}(f^2-\E f^2)(X_i))>\\
&>\bar{A}(N^{-1}\gamma^2_2(\ccF)+N^{-\frac{1}{2}}\alpha\gamma_2(\ccF))+\bar{B}N^{-\frac{1}{2}}\alpha^2 u)\ls \exp(-\min(\sqrt{N u},u^2)),
\end{align*}
where $\bar{A},\bar{B}$ are universal constants.
\end{theo}
\begin{dwd}
We mimic the proof of Proposition \ref{pro1}.  Assume there is an admissible sequence of nets $(\ccF_n)^{\infty}_{n=0}$
such that $\ccF_n\subset \ccF_{n+1}$ and $|\ccF_n|\ls N_n=\varphi_2(2^n)=2^{2^{2n}}-1$ for $n\gs 1$ and $|\ccF_0|=1$
and $\bigcup^{\infty}_{n=0}\ccF_n$ is dense in $\ccF$.
There exists $k_0\gs 0$ such that $2^{2k_0}\ls N< 2^{2(k_0+1)}$ and hence
$$
\varphi^{-1}(N_n)= \left\{\begin{array}{ll}
\frac{1}{\sqrt{N}}2^{2n} & n> k_0 \\
2^{n} & n\ls k_0
\end{array} \right.
$$
Therefore $\varphi^{-1}(N_n)$ behaves like $\varphi^{-1}_2(N_n)$ for $n\ls k_0$ and $N^{-\frac{1}{2}}\varphi^{-1}_1(N_n)$ for $n> k_0$. Fix $f\in \ccF$. The main idea is to apply the following split
\begin{align*}
&f^2(X_i)=(f-\pi_{k_0}(f)+\pi_{k_0}(f))^2(X_i)=(f(X_i)-\pi_{k_0}(f))^2+\\
&+2(f-\pi_{k_0}(f))\pi_{k_0}(f)(X_i)+(\pi_{k_0}(f))^2(X_i).
\end{align*}
which implies that
\begin{align*}
&\frac{1}{N}\sum^N_{i=1}f^2(X_i)\ls \frac{1}{N}\sum^N_{i=1}(f-\pi_{k_0}(f))^2(X_i)+\\
&+2\frac{1}{N}\sum^N_{i=1}(f-\pi_{k_0}(f))\pi_{k_0}(f)(X_i)+\frac{1}{N}\sum^N_{i=1}(\pi_{k_0}(f))^2(X_i).
\end{align*}
Let us define
\begin{align*}
& P_N(f)=\frac{1}{N}\sum^N_{i=1}(f-\pi_{k_0}(f))^2(X_i),\\
& Q_N(f)=\frac{2}{N}\sum^{N}_{i=1}((f-\pi_{k_0}(f))\pi_{k_0}(f))(X_i), \\
& R_N(f)=\frac{1}{N}\sum^N_{i=1}(\pi_{k_0}(f))^2(X_i).
\end{align*}
Consequently
$$
S_N(f)= P_N(f)+Q_N(f)+R_N(f)-\E P_N(f) -\E Q_N(f) -\E R_N(f). 
$$
Now we aim to study $P_N(f),Q_N(f)$ and $R_N(f)$ separately.
\smallskip
 
\noindent
We start from $P_N(f)$. The main point is to explore our concept of certain structure $n_k$, $k\gs 0$ that
may depend on $f$ to suitably bound $P_N(f)$. We use the following chaining  
$$
(f-\pi_{k_0}(f))(X_i)= \sum^{\infty}_{k=0}(\pi_{n_k}(f)-\pi_{n_{k+1}}(f))(X_i),
$$
where $(n_k)^{\infty}_{k=0}$ is defined as follows $n_0(f)=k_0$ and
$$
n_i= \inf\{n>n_{i-1}:\;2d(f,\ccF_n)< d(f,\ccF_{n_{i-1}})\}.
$$ 
Let $\Delta_k(f)=\frac{1}{N}\sum^N_{i=1}(\pi_{n_k}(f)-\pi_{n_{k+1}}(f))^2(X_i)$.
By the triangle inequality it is clear that
\be\label{raus}
[\frac{1}{N}\sum^N_{i=1}(f-\pi_{k_0}(f))^2(X_i)]^{\frac{1}{2}}\ls \sum^{\infty}_{k=0} [\Delta_k(f)]^{\frac{1}{2}}. 
\ee
Moreover
$$
[\Delta_k(f)]^{\frac{1}{2}}\ls [(\Delta_k(f) -\E \Delta_k(f))_{+}]^{\frac{1}{2}}+[\E \Delta_k(f)]^{\frac{1}{2}}.
$$
Clearly 
\be\label{raus2}
[\E \Delta_k(f)]^{\frac{1}{2}}\ls Cd(\pi_{n_{k}}(f),\pi_{n_{k+1}}(f)).
\ee
By the Bernstein inequality
$$
\|(\Delta_k(f) -\E \Delta_k(f))_{+}\|_{\varphi}\ls  
K^2N^{-\frac{1}{2}}d(\pi_{n_{k}}(f),\pi_{n_{k+1}}(f))^2.
$$
Therefore
\begin{eqnarray}
&& [(\Delta_k(f) -\E \Delta_k(f))_{+}]^{\frac{1}{2}}\ls \nonumber \\
&&\ls KN^{-\frac{1}{4}}d(\pi_{n_{k}}(f),\pi_{n_{k+1}}(f))(3\varphi^{-1}(N_{k+1})+\varphi^{-1}(Z_1))^{\frac{1}{2}}\ls \nonumber\\
\label{burak1}&&\ls KN^{-\frac{1}{4}}d(\pi_{n_{k}}(f),\pi_{n_{k+1}}(f))(2[\varphi^{-1}(N_{k+1})]^{\frac{1}{2}}+[\varphi^{-1}(Z_1)]^{\frac{1}{2}}),
\end{eqnarray}
where $Z_1$ is given by
$$
Z_1=\sum^{\infty}_{k=1} \frac{1}{N_k^3}U_k,\;\;U_k=\sum_{g,h\in \ccF_k} \varphi(X(g,h))
$$
and
$$
X(g,h)=\frac{\sum^N_{i=1} [(g-h)^2-\E(g-h)](X_i)}{K^2 N^{\frac{1}{2}}d(g,h)^2}.
$$
Clearly $\E Z_1\ls 1$ by (\ref{eq:5.5}).
We can sum up the bounds, i.e. using (\ref{raus}), (\ref{raus2})  and (\ref{burak1}) we get
\begin{align*}
& [P_N(f)]^{\frac{1}{2}}\ls \sum^{\infty}_{k=0} d(\pi_{n_{k}}(f),\pi_{n_{k+1}}(f))(C+2KN^{-\frac{1}{4}}[\varphi^{-1}(N_{k+1})]^{\frac{1}{2}})+\\
&+ KN^{-\frac{1}{4}}\sum^{\infty}_{k=0} d(\pi_{n_{k}}(f),\pi_{n_{k+1}}(f))[\varphi^{-1}(Z_1)]^{\frac{1}{2}}).
\end{align*}
We have noticed that $\varphi^{-1}(N_n)=N^{-\frac{1}{2}}2^{2n}$, and by the construction
$$
d(\pi_{n_{k+1}}(f),\pi_{n_{k}}(f))\ls d(f,\pi_{n_k}(f))+2d(f,\pi_{n_{k}-1}(f))
$$
and for $k\gs 0$
$$
d(\pi_{n_{k}}(f),f)\gs 2 d(\pi_{n_{k+1}}(f),f).
$$
It implies that
\begin{align*}
& [P_N(f)]^{\frac{1}{2}}\ls 12KN^{-\frac{1}{2}}\sigma_{k_0}(f)+6Cd(f,\pi_{k_0}(f))+\\
&+6CN^{-\frac{1}{4}}d(f,\pi_{k_0}(f))\varphi^{-1}(Z_1)^{\frac{1}{2}},
\end{align*}
where $\sigma_{k_0}(f)=\sum^{\infty}_{n=k_0}2^nd(t,\ccF_n)$. Obviously
$$
\sup_{f\in \ccF}\sigma_{k_0}(f)\ls \gamma_2(\ccF).
$$
Moreover
\be\label{eq:5.6}
d(f,\pi_{k_0}(f))\ls d(f,\pi_{k_0}(f))(2^{k_0+1}N^{-\frac{1}{2}})\ls 2N^{-\frac{1}{2}}\gamma_2(\ccF)
\ee
and on the other hand
\be\label{eq:5.7}
d(f,\pi_{k_0}(f))\ls 2\alpha.
\ee
Since $\bigcup^{\infty}_{n=0}\ccF_n$ is dense in $\ccF$ we acquire
for each $f\in \ccF$
\begin{eqnarray}
&& P_N(f)\ls [12(K+C)N^{-\frac{1}{2}}\gamma_2(\ccF)+\nonumber \\
&& +12CN^{-\frac{1}{4}}\min\{\alpha,N^{-\frac{1}{4}}\gamma_2(\ccF)\}[\varphi^{-1}(Z_1)]^{\frac{1}{2}}]^2\ls 2(12)^2(K+C)^2 N^{-1}\gamma_2^2(\ccF)+\nonumber\\
\label{bleble} &&+ 2(12)^2 C^2 N^{-\frac{1}{2}}\min\{\alpha^2,N^{-\frac{1}{2}}\gamma^2_2(\ccF)\}
\varphi^{-1}(Z_1).
\end{eqnarray}
Moreover by (\ref{eq:5.6}) and (\ref{eq:5.7})
\be\label{eq:5.8}
\E P_N(f)\ls C^2d(f,\pi_{k_0}(f))^2\ls 4C^2\min\{\alpha^2,N^{-1}\gamma^2_2(\ccF)\}.
\ee
The second point is to consider $R_N(f)$. We use the following chaining
\begin{align*}
& \frac{1}{N}\sum^N_{i=1} ((\pi_{k_0}(f))^2-\E (\pi_{k_0}(f))^2-(\pi_0(f))^2+\E (\pi_0(f))^2)(X_i) \ls \\
&\ls \sum^{\infty}_{k=1} \frac{1}{N} \sum^N_{i=1} ((\pi_{n_k}(f))^2-(\pi_{n_{k-1}}(f))^2-\E (\pi_{n_k}(f))^2+\E (\pi_{n_{k-1}}(f))^2)(X_i).
\end{align*}
where
$$
n_i(f)=\min\{\inf\{n>n_{i-1}(f):\; 2d(f,\ccF_n)<d(f,\ccF_{n_{i-1}})\},k_0\}.
$$
Note that $n_i=k_0$ for a finite $i\gs 0$. Using (\ref{eq:5.6}) we get
\begin{align*}
&| \frac{1}{N} \sum^N_{i=1} ((\pi_{n_k}(f))^2-(\pi_{n_{k-1}}(f))^2-\E (\pi_{n_k}(f))^2+\E (\pi_{n_{k-1}}(f))^2)(X_i)|\ls \\
&\ls 2\alpha N^{-\frac{1}{2}}K^2 d(\pi_{n_k}(f),\pi_{n_{k-1}}(f))(3\varphi^{-1}(N_{n_k})+\varphi^{-1}(Z_2)),
\end{align*}
where 
$$
Z_2=\sum^{\infty}_{k=1} \frac{1}{N^3_k}V_k,\;\;V_k=\sum_{g,h\in \ccF_k}
\varphi(X(g,h)),
$$
and
$$
X(g,h)=\frac{\sum^N_{i=1}(g^2-h^2-\E g^2+\E h^2)}{2K^2N^{\frac{1}{2}}\alpha d(g,h)}.
$$
Note that $\E Z_2$ due to (\ref{burak3}). By the construction of $n_k$, $k\gs 1$
$$
d(\pi_{n_k}(f),\pi_{n_{k-1}}(f))\ls d(f,\pi_{n_k}(f))+2d(f,\pi_{n_{k}-1}(f))
$$
and
$$
d(f,\pi_{n_k}(f))\gs 2d(f,\pi_{n_{k-1}}(f)).
$$
Since $\varphi^{-1}(N_n)=2^n$ we obtain that
\begin{align*}
& \frac{1}{N}\sum^N_{i=1} ((\pi_{k_0}(f))^2-\E (\pi_{k_0}(f))^2-(\pi_0(f))^2+\E (\pi_0(f))^2)(X_i) \ls \\
&\ls 12 K^2N^{-\frac{1}{2}} \alpha \sigma_{0,k_0}(f) +6N^{-\frac{1}{2}}K^2\alpha d(f,\pi_0(f))\varphi^{-1}(Z_2),
\end{align*}
where 
$$
\sigma_{0,k_0}(f)=\sum^{k_0}_{n=0}2^n d(f,\ccF_n)\ls \gamma_2(\ccF). 
$$
Obviously $\E Z_2\ls 1$.  On the other hand 
$$
\frac{1}{N}\sum^N_{i=1} ((\pi_0(f))^2-\E (\pi_0(f))^2)(X_i)\ls K^2 N^{-\frac{1}{2}}\alpha^2\varphi^{-1}(Z_3)
$$
where
$$
Z_3=\varphi(\frac{\sum^N_{i=1}((\pi_0(f))^2-\E(\pi_0(f))^2)(X_i)}{K^2 N^{\frac{1}{2}}\alpha^2}).
$$
Again observe that $\E Z_3$ by (\ref{burak4}).
It implies that
\be\label{burak5}
R_N(f)-\E R_N(f)\ls 12K^2 N^{-\frac{1}{2}}\alpha \gamma_2(\ccF)+6K^2 N^{-\frac{1}{2}}\alpha^2(\varphi^{-1}(Z_2)+\varphi^{-1}(Z_3)).
\ee
Moreover
\be\label{rain}
\E R_N(f)\ls C^2 d^2(\pi_{k_0}(f),f)\ls 4C^2 \min\{\alpha^2,N^{-1}\gamma^2_2(\ccF)\}. 
\ee
It remains to bound $Q_N(f)$. Clearly
\be\label{rain1}
Q_N(f)\ls 2[P_N(f)]^{\frac{1}{2}}[R_N(f)]^{\frac{1}{2}}.
\ee
and
$$
[R_N(f)]^{\frac{1}{2}}\ls [(R_N(f)-\E R_N(f))_{+}]^{\frac{1}{2}}+[\E R_N(f)]^{\frac{1}{2}}.  
$$
Therefore using (\ref{rain}) and (\ref{rain1}) 
\begin{align*}
& Q_N(f)\ls 2([\E R_N(f)]^{\frac{1}{2}}+[(R_N(f)-\E R_N(f))_{+}]^{\frac{1}{2}})
[P_N(f)]^{\frac{1}{2}} \ls \\
& \ls 2(2C\min\{\alpha,N^{-\frac{1}{2}}\gamma_2(\ccF)\} +
[(R_N(f)-\E R_N(f))_{+}]^{\frac{1}{2}})[P_N(f)]^{\frac{1}{2}}.
\end{align*} 
and consequently  
\begin{eqnarray}
&& Q_N(f)\ls 4C\min\{\alpha,N^{-\frac{1}{2}}\gamma_2(\ccF)\}[P_N(f)]^{\frac{1}{2}}+\nonumber\\
&&+2(R_N(f)-\E R_N(f) +P_N(f)).
\end{eqnarray}
Moreover by (\ref{rain}) and (\ref{eq:5.8})
\be\label{burak7}
|\E Q_N(f)|\ls 2[\E P_N(f)]^{\frac{1}{2}}[\E R_N(f)]^{\frac{1}{2}}\ls 8C^2 \min\{\alpha^2,N^{-1}\gamma^2_2(\ccF)\}.
\ee
Finally we can sum up the inequalities. First observe that  
\begin{align*}
& P_N(f)+Q_N(f)+R_N(f)-\E(P_N(f)+Q_N(f)+R_N(f))\ls PN(f)+R_N(f)- \\
&\ls -\E R_N(f)+2([(R_N(f)-\E R_N(f))_{+}]^{\frac{1}{2}}+[\E R_N(f)]^{\frac{1}{2}}) [P_N(f)]^{\frac{1}{2}}+|\E Q_N(f)|\ls\\
&\ls 2P_N(f)+2(R_N(f)-\E R_N(f))+2[\E R_N(f)]^{\frac{1}{2}} [P_N(f)]^{\frac{1}{2}}+|\E Q_N(f)|
\end{align*}
and hence using (\ref{bleble}), (\ref{burak5}) and (\ref{burak7})
\begin{align*}
& P_N(f)+Q_N(f)+R_N(f)-\E(P_N(f)+Q_N(f)+R_N(f))\ls \\
&\ls 2[12(K+C)N^{-\frac{1}{2}}\gamma_2(\ccF)+12CN^{-\frac{1}{4}}\min\{\alpha,N^{-\frac{1}{4}}\gamma_2(\ccF)\}[\varphi^{-1}(Z_1)]^{\frac{1}{2}}]^2+\\
&+ 2[12K^2 N^{-\frac{1}{2}}\alpha \gamma_2(\ccF)+K^2 N^{-\frac{1}{2}}\alpha^2(6\varphi^{-1}(Z_2)+2\varphi^{-1}(Z_3))]+\\
&+4C\min\{\alpha,N^{-\frac{1}{2}}\gamma_2(\ccF)\}
[12(K+C)N^{-\frac{1}{2}}\gamma_2(\ccF)+\\
&+12CN^{-\frac{1}{4}}\min\{\alpha,N^{-\frac{1}{4}}\gamma_2(\ccF)\}[\varphi^{-1}(Z_1)]^{\frac{1}{2}}]+\\
&+8C^2\min\{\alpha^2,N^{-1}\gamma^2_2(\ccF)\}.
\end{align*}
Therefore
\begin{align*}
& P_N(f)+Q_N(f)+R_N(f)-\E(P_N(f)+Q_N(f)+R_N(f))\ls \\
&\ls (24)^2 (K+C)^2 N^{-1}\gamma^2_2(\ccF)+24C(K+C)\min\{\alpha,N^{-\frac{1}{2}}\gamma_2(\ccF)\}N^{-\frac{1}{2}}\gamma_2(\ccF)+\\
&+24K^2N^{-\frac{1}{2}}\alpha \gamma_2(\ccF)+8C^2\min\{\alpha^2,N^{-1}\gamma^2_2(\ccF)\}+\\
&+(24)^2C^2N^{-\frac{1}{2}}\min\{\alpha^2,N^{-\frac{1}{2}}\gamma^2_2(\ccF)\}\varphi^{-1}(Z_1)+\\
&+6K^2N^{-\frac{1}{2}}\alpha^2(\varphi^{-1}(Z_2)+\varphi^{-1}(Z_3)).
\end{align*}
It proves that
\begin{align*}
& P_N(f)+Q_N(f)+R_N(f)-\E(P_N(f)+Q_N(f)+R_N(f))\ls \\
&\ls  A(N^{-1}\gamma^2_2(\ccF)+N^{-\frac{1}{2}}\alpha \gamma_2(\ccF))+BN^{-\frac{1}{2}}\varphi^{-1}(Z),
\end{align*}
where $A=(24)^2K^2+24C(K+C)+8C^2\ls (30)^2(K+C)$ and $B=6(24)^2\max\{C^2,K^2\}$ and
$$
Z=\frac{1}{2}Z_1+\frac{1}{4}Z_2+\frac{1}{4}Z_3.
$$
Not that using that $D^{-1}\varphi(x)\gs \varphi(D^{-1}x)$ for $D\gs 1$
\begin{align*}
 &\frac{1}{2}\varphi^{-1}(Z_1)+\frac{1}{4}\varphi^{-1}(Z_2)+\frac{1}{4}\varphi^{-1}(Z_3)\ls
\varphi^{-1}(\frac{1}{2}Z_1)+\varphi^{-1}(\frac{1}{4}Z_2)+\varphi^{-1}(\frac{1}{4}Z_3)\ls\\
&\ls 3\varphi^{-1}(\frac{1}{2}Z_1+\frac{1}{4}Z_2+\frac{1}{4}Z_3)=3\varphi^{-1}(Z).
\end{align*}
Clearly $\E Z\ls 1$ which ends the proof.
\end{dwd}
 
\section{Compressed sensing}

In the compressed sensing we consider $N\times M$  matrix $A$, where $N<<M$.
We want to reconstruct all vectors
$x\in \R^M$ of sparse support, i.e. all vectors such that 
$|\{i\in \{1,...,M\}:x_i\neq 0\}|\ls m$. For simplicity denote
by $\Sigma_m$ the space of $m$-sparse vectors in $\R^m$.
The main tool to recover any sparse vector $x$ is the $l^1$
minimization which works \cite{Can} whenever $A$ satisfies Restricted Isometry Property, namely
$$
1-\delta_m\ls |Ax|^2_2\ls 1+\delta_m
$$
for all $m$-sparse $x$ in the unit sphere of $\R^M$ (which we denote by $S_2(\Sigma_m)$). 
The main result for RIP is the following
\bt
Let $1\ls m\ls N/2$
Let $A$ be an $N\times M$ matrix. If
$\delta_{2m}(A)<\sqrt{2}-1$, then $A$ satisfies the exact reconstruction property of order $m$
(all vectors that are $m$-sparse are reconstructed).
\et
Consider $A=N^{-\frac{1}{2}}(Y_1,...,Y_N)^T$ and let $Y_i$ be independent subgaussians with 
$\|Y_i\|_{\varphi_2}\ls \alpha$. 
It is clear then
$$
(1-\delta_{2m})|x|^2_2\ls |Ax|^2_2\ls (1+\delta_{2m})|x|^2_2, 
$$
where
$$
\delta_{2m}=\sup_{x\in S_2(\Sigma_{2m})}
|\frac{1}{N}\sum^N_{i=1}\langle Y_i,x\rangle^2-1|.
$$
To ensure $\delta_{2m}<\sqrt{2}-1$ we need that 
$$
\sup_{x\in S_2(\Sigma_{2m})}
|\frac{1}{N}\sum^N_{i=1}\langle Y_i,x\rangle^2-1|<\sqrt{2}-1.
$$
Note that for any $x,y\in S_2$ by the Schwarz inequality
$$
\| \langle Y_i, x-y\rangle \|_{\varphi_2}\ls \|Y_i\|_{\varphi_2}\|x-y\|_2\ls \alpha\|x-y\|_2.
$$
Let $\ccF=\{f(\cdot)=\langle x,\cdot\rangle:\;x\in S_2(\Sigma_{2m})\}$
Therefore for $f(\cdot)=\langle x,\cdot\rangle$ and $g(\cdot)=\langle y ,\cdot\rangle$
\be\label{bam}
d(f,g)\ls \alpha\|x-y\|_2, \;\;f,g\in \ccF.
\ee
We can use the result of the previous section
provided we can compute $\gamma_2(S_2(\Sigma_{2m}))$.
By the majorizing 
measure theorem 
$$
\gamma_2(S_2(\Sigma_{2m}))\sim \E(\sup_{t\in
S_2(\Sigma_{2m})}\sum^M_{i=1}t_ig_i,
$$ 
where $g_i$, $1\ls i\ls M$ are independent standard Gaussian variables.
Moreover
$$
\E(\sup_{t\in
S_2(\Sigma_{2m})}\sum^M_{i=1}t_ig_i= \E(\sum^{2m}_{i=1}(g_i^{\ast})^2)^{\frac{1}{2}},
$$ 
where $g_i^{\ast}$ is a non-decreasing rearrangement of $(|g_i|)^{N}_{i=1}$.
Finally we have a simple result to compute the last quantity.
\bl
Let
There exists absolute constants $c_0,c_1>0$ such that
the following holds. Let $(g_i)^M_{i=1}$ be a family of independent standard Gaussian 
random variables and $(g_i^{\ast})^M_{i=1}$ its non-decreasing rearrangement. For any
$k\ls M$ we have
$$
\sqrt{c_1k\log(k^{-1}M)}\ls
\E(\sum^{k}_{i=1}(g_i^{\ast})^2)^{\frac{1}{2}}\ls 2\sqrt{k\log(c_0k^{-1}M)}.
$$
\el
It shows that 
$$
\gamma_2(S_2(\Sigma_{2m}))\ls K_0 \sqrt{m\log(c_0 m^{-1}M)},
$$
where $K_0,c_0$ are absolute constants. Therefore due to (\ref{bam})
$$
\gamma_2(\ccF)\ls K_0\alpha\sqrt{m\log(c_0m^{-1}M)}.
$$
Denote $\gamma=K_0(m\log(c_0m^{-1}M))^{\frac{1}{2}}$, by Theorem \ref{thm1} we get
\begin{align*}
& \P(\sup_{x\in S_2(\Sigma_{2m})}
|\frac{1}{N}\sum^N_{i=1}\langle Y_i,x\rangle^2-1|>
A(N^{-1}\alpha^2\gamma^2+N^{-\frac{1}{2}}\alpha^2\gamma)+
BN^{-\frac{1}{2}}\alpha^2u)\ls\\
&\ls \exp(-(\log 2)\min\{N^{-\frac{1}{2}}u,u^2\}),
\end{align*}
where $1\ls A,B<\infty$. Now it is clear that we should take $u=c_0\gamma$ and then
for any 
$$
(A+B+c_0)N^{-\frac{1}{2}}\gamma\alpha^2=\delta\ls 1,
$$
which in particular means that  $\gamma\ls N^{-\frac{1}{2}}$ the following inequality holds
$$
\P(\sup_{x\in S_2(\Sigma_{2m})}
|\frac{1}{N}\sum^N_{i=1}\langle Y_i,x\rangle^2-1|>\delta)\ls  \exp(-(\log 2)\gamma^2).
$$


\begin{thebibliography}{99}

\footnotesize


\bibitem{Can} \textsc{Candes, E.J.} (2008), The restricted isometry property and its implications for compressed sensing. \textit{C.R.Math.Acad.Sci.Paris}, {\bf 346}, N.9-10, 589-592.
\bibitem{MPT} \textsc{Mendelson, S. Pajor, A and Tomczak-Jaegermann} (2008) Reconstruction and subgaussian operators in asymptotic geometric analysis, \textit{Geom. Funct. Anal.} {\bf 17}
N. 4, 1248-1282.
\bibitem{PaM} \textsc{Mendelson, S, Pouris, G} (2011) On generic chaining and the smallest singular value of random matrices with heavy tails. Preprint.
\bibitem{Tal} \textsc{Talagrand, M.} (2005), The generic chaining. \textit{Springer-Verlag}.

\end{thebibliography}
\end{document}